\documentclass[11pt]{article}
\usepackage{amsfonts}
\usepackage{latexsym,amssymb,amsmath,graphics,cite,arydshln}
\usepackage{tikz}
\usepackage{fancyhdr}
\usepackage{lipsum}
\usepackage{lineno}
\usepackage{color}

\usepackage{algorithm} 
\usepackage{algorithmic} 
\usepackage{multirow} 
\usepackage{amsmath}
\usepackage{xcolor}

\begin{document}
\newcommand{\qed}{\hphantom{.}\hfill $\Box$\medbreak}
\newcommand{\proof}{\noindent{\bf Proof \ }}

\newtheorem{theorem}{Theorem}[section]
\newtheorem{lemma}[theorem]{Lemma}
\newtheorem{corollary}[theorem]{Corollary}
\newtheorem{remark}[theorem]{Remark}
\newtheorem{example}[theorem]{Example}
\newtheorem{definition}[theorem]{Definition}
\newtheorem{construction}[theorem]{Construction}
\newtheorem{fact}[theorem]{Fact}
\newtheorem{proposition}[theorem]{Proposition}
\newtheorem{conjecture}[theorem]{Conjecture}


	\title{New results of Bollob\'{a}s-type theorem for affine subspaces and projective subspaces}
	
	\author{\small Shuhui\ Yu, Xin\ Wang \footnote{S.Yu is with the Department of Mathematics and Physics, Suzhou Vocational University, Suzhou 215000, P. R. China  (e-mail: yushuhui\_suda@163.com); X. Wang  is with the Department of Mathematics, Soochow University, Suzhou 215006, P. R. China  (e-mail: xinw@suda.edu.cn)}}
	\date{}
	\maketitle
	\begin{abstract}

Bollob\'{a}s-type theorem has received a lot of attention due to its application in graph theory. In 2015, G\'{a}bor Heged{\"u}s gave an upper bound of bollob\'{a}s-type affine subspace families for $q\neq 2$, and constructed an almost sharp affine subspaces pair families. In this note, we prove a new upper bound for bollob\'{a}s-type affine subspaces without the requirement of $q\neq 2$, and construct a pair of families of affine subspaces, which shows that our upper bound is sharp. We also give an upper bound for  bollob\'{a}s-type projective subspaces, and prove that the Heged{\"u}s's conjecture holds when $q=2$.

\smallskip
	\end{abstract}
	
	\section{Introduction}
	
	Bollob\'{a}s proved in \cite{Bo1965} the following famous result.

\begin{theorem}{\rm\cite{Bo1965}}\label{Bo1}
	Let $A_{1},A_{2},\ldots, A_{m}$ and $B_{1},B_{2},\ldots,B_{m}$ be two families of sets such that $A_{i}\cap B_{j}=\emptyset$ if and only if $i=j$. Then $$\sum_{i=1}^{m}\frac{1}{{
 |A_{i}+|B_{i}| \choose
 |A_{i}| }}\leq 1.$$
In particular, if $|A_{i}|=r$ and $|B_{i}|=s$ for each $1\leq i\leq m$, then $$m\leq {r+s\choose s}.$$
	\end{theorem}

The following uniform version of Bollob\'{a}s's theorem was proved by Lov\'{a}sz in \cite{Lo1977} using tensor product method.

\begin{theorem}{\rm\cite{Lo1977}}\label{Bo2}
	Let $\mathcal{F}=\{A_{1},\ldots,A_{m}\}$ be an $r$-uniform family and $\mathcal{G}=\{B_{1},\ldots,B_{m}\}$ be an $s$-uniform family such that $A_{i}\cap B_{i}=\emptyset$ for each $1\leq i\leq m$ and $ A_{i}\cap B_{j}\neq\emptyset$ for any $1\leq i<j\leq m$ Then  $$m\leq \begin{pmatrix}
 r+s\\
r \\
\end{pmatrix}.$$
	\end{theorem}

Lov\'{a}sz also proved the following generalization of Bollob\'{a}s'theorem for subspaces in \cite{Lo1979}.

\begin{theorem}{\rm\cite{Lo1979}}\label{Bo3}
	Let $\mathbb{F}$ be an arbitrary field and $V$ be an $n$-dimensional vector space over the field $\mathbb{F}$. Let $U_{1},U_{2},\ldots,U_{m}$ denote $r$-dimensional subspaces of $V$ and $V_{1},V_{2},\ldots,V_{m}$ denote $s$-dimensional subspaces of the vector space $V$. Assume that $U_{i}\cap V_{i}=\{0\}$ for each $1\leq i\leq m$ and $U_{i}\cap V_{j}\neq \{0\}$ whenever $1\leq i<j\leq m$. Then $$m\leq \begin{pmatrix}
 r+s\\
r \\
\end{pmatrix}.$$
	\end{theorem}

	Heged{\"u}s gave an affine version of Theorem \ref{Bo3} in \cite{He2015} by linear algebra method.

\begin{theorem}{\rm\cite{He2015}}\label{Bo4}
	Let $A_{1},A_{2},\ldots,A_{m}$ and $B_{1},B_{2},\ldots,B_{m}$ be affine subspaces of an $n$-dimensional affine space $W$ over the finite field $\mathbb{F}_{q}$, where $q\neq 2$. Assume that $\{(A_{i},B_{i})\colon 1\leq i\leq m\}$ is cross-intersecting (i.e., $A_{i}\cap B_{i}=\emptyset$ for each $1\leq i \leq m$ and $A_{i}\cap B_{j}\neq \emptyset$ whenever $1\leq i<j\leq m)$. Then $m\leq q^{n}+1$. 
	\end{theorem}

In \cite{He2015}, let $m(n,q)$ denote the maximal size of a cross-intersecting pair of families of affine subspaces $\{(A_{i},B_{i})\}_{1\leq i\leq m}$. Heged{\"u}s also constructed a cross-intersecting pair of families of affine subspace, leading to a lower bound of $m(n,q)$.

\begin{theorem}{\rm\cite{He2015}}\label{Bo5}
	Let $n\geq 1$ and $q$ be an arbitrary prime power. Then $m(n,q)\geq \frac{q^{n}-1}{q-1}$.
	\end{theorem}

In \cite{He2015}, Heged{\"u}s posed two interesting open problems:

Open problem 1. What can we say about $m(n,2)$?

Open problem 2.  What is the precise value of $m(n,q)$, if $q>2$?

Heged{\"u}s also conjectured the following projective version of Theorem \ref{Bo4}.

\begin{conjecture}\label{HegedusConj}
Let $\mathbb{F}$ be an arbitrary field. Let $A_{1},\ldots,A_{m}$ and $B_{1},\ldots,B_{m}$ be projective subspaces of an $n$-dimensional projective space $W$ over the field $\mathbb{F}$. Assume
that $\{(A_{i},B_{i})\colon 1\leq i\leq m\}$ is cross-intersecting $($i.e. $A_{i}\cap B_{i}=\emptyset$ for each $1\leq i\leq m$ and $A_{i}\cap B_{j}\neq \emptyset$ whenever $1\leq i<j\leq m)$. Then
$m\leq 2^{n+1}-2.$
\end{conjecture}

In 2024, Heged{\"u}s and Frankl proved that $m\leq 2^{n+1}$ by using tensor product method \cite{HF2024}.

	The first part of this note is devoted to obtain a new upper bound of $m(n,q)$, and construct a perfect cross-intersecting pair of affine subspaces. It is proved that $m(n,q)= 2\frac{q^{n}-1}{q-1}$ for any prime power $q$ and positive integer $n$. Our result not only improve Theorem \ref{Bo5}, but also reduce the requirement of $q\neq 2$, which solved the open problems presented by  Heged{\"u}s. The second part of this note is devoted to obtain an upper bound of cross-intersecting pair of projective subspaces, It is proved that the conjecture of Heged{\"u}s is hold for $q=2$.

	\section{New upper bound of $m(n,q)$}
	
	\begin{theorem}\label{Main}
	Let $A_{1},A_{2},\ldots,A_{m}$ and $B_{1},B_{2},\ldots,B_{m}$ be affine subspaces of an affine space $W$ with dimension $n$ over the finite field $\mathbb{F}_{q}$. Assume that $\{(A_{i},B_{i})\colon 1\leq i\leq m\}$ is cross-intersecting $($i.e. $A_{i}\cap B_{i}=\emptyset$ for each $1\leq i\leq m$ and $A_{i}\cap B_{j}\neq \emptyset$ whenever $1\leq i<j\leq m)$. Then $m\leq 2\cdot \frac{q^{n}-1}{q-1}$. 
	\end{theorem}
	
\begin{proof}
Assume that $A_{i}=v_{i}+V_{i}$ and $B_{i}=u_{i}+U_{i}$ for any $1\leq i\leq m$, where $V_{i}$ and $U_{i}$ are subspaces of $W$, $v_i$ and $u_i$ are some vectors from $W$. Since $A_{i}\cap B_{i}=\emptyset$, we know that $v_{i}-u_{i}\notin V_{i}-U_{i}$ and $V_{i}-U_{i}\neq W$. Then there exists a hyperplane $H_{i}$ of $W$ such that $v_{i}-u_{i}\notin H_{i}$ and  $V_{i}-U_{i}\subseteq H_{i}$. Notice that replacing $v_{i}+V_{i}$ and $u_{i}+U_{i}$ by $v_{i}+H_{i}$ and $u_{i}+H_{i}$ respectively does not violate the cross-intersecting property. Without loss of generality, suppose that $A_{i}=v_{i}+H_{i}$ and $B_{i}=u_{i}+H_{i}$ for any $1\leq i\leq m$, where $H_{i}$ is a hyperplane.

Claim: At most two distinct affine subspaces (cosets) of the same hyperplane $H$ can appear in $\{A_{1},A_{2},\ldots,A_{m}\}.$

Suppose for contradiction that there exist three distinct affine subspaces $v_{1}+H,v_{2}+H,v_{3}+H$ with the same hyperplane $H$ appear in $\{A_{1},A_{2},\ldots,A_{m}\}.$ Without loss of generality, we assume that $v_{1}+H=A_{i}$, $v_{2}+H=A_{j}$ and $v_{3}+H=A_{\ell}$, where $1\leq i<j<\ell\leq m$. Since $\{(A_{i},B_{i})\}_{1\leq i\leq m}$ is cross-intersecting and $A_{i}\cap B_{i}=\emptyset$, there exist three affine subspaces $v_{4}+H$, $v_{5}+H$ and $v_{6}+H$  such that $v_{4}+H=B_{i}$, $v_{5}+H=B_{j}$ and  $v_{6}+H=B_{\ell}$.  Since $A_{i}\cap B_{j}\neq \emptyset$ whenever $1\leq i<j\leq m$, we have $v_{1}-v_{6}\in H$ and $v_{2}-v_{6}\in H$. It follows that $v_{1}+H=v_{6}+H$ and $v_{2}+H=v_{6}+H$, which is a contradiction because $v_{1}+H$ and $v_{2}+H$ are different cosets of $H$.

Since the number of hyperplanes of $\mathbb{F}_{q}^{n}$ is $\frac{q^{n}-1}{q-1}$, by Claim, we have $m\leq 2\cdot \frac{q^{n}-1}{q-1}.$\qed

\end{proof}

In the following, we generalize the construction of cross-intersecting pair of families of affine subspaces in \cite{He2015} to give an exact value of $m(n,q)$.

\begin{theorem}\label{Main2}
	Let $n\geq 1$ and $q$ be an arbitrary prime power. Then $m(n,q)= 2\cdot \frac{q^{n}-1}{q-1}$.
	\end{theorem}
\begin{proof}
Let $t=\frac{q^{n}-1}{q-1}$, $\mathcal{H}=\{H_{1},H_{2},\ldots,H_{t}\}$ be  the set of all hyperplanes of the vector space $\mathbb{F}_{q}^{n}$. For each $1\leq i\leq t$, we fix a vector $\beta_{i}\in \mathbb{F}_{q}^n\backslash H_{i}$. Define $$A_{i}:=H_{i},\ A_{t+i}:=H_{i}+\beta_{i},$$ and $$B_{i}:=H_{i}+\beta_{i},\ B_{t+i}:=H_{i}.$$

Clearly $A_{i},A_{t+i},B_{i},B_{t+i}$ are affine subspaces  for each $1\leq i\leq t$. Since $\beta_{i}\not\in H_{i}$ for each $1\leq i\leq t$, we have that $A_{i}\cap B_{i}=\emptyset$ and $A_{m+i}\cap B_{t+i}=\emptyset$.

On the other hand, since $\beta_{i}\in H_{i}-H_{j}=\mathbb{F}_{q}^{n}$ and $\dim(H_{i}\cap H_{j})\geq 2(n-1)-n= n-2$, it follows that $A_{i}\cap B_{j}\neq \emptyset$ for each $1\leq i<j\leq t$, $A_{i}\cap B_{t+j}\neq\emptyset$ for each $1\leq i,j\leq t$, and $A_{t+i}\cap B_{t+j}\neq \emptyset$ for each $1\leq i<j\leq t.$ 

Hence, $\{(A_i,B_i)\colon 1\leq i\leq t\}\cup \{(A_{t+i},B_{t+i})\colon 1\leq i\leq t\}$ is cross-intersecting. By Theorem \ref{Main} we have $m= 2\cdot \frac{q^{n}-1}{q-1}$. \qed
\end{proof}

\section{Cross-intersecting of projective subspaces}

In the following, we will prove that Heged{\"u}s' conjecture holds for $q=2$.

\begin{theorem}\label{Main2}
Let $A_{1},\ldots,A_{m}$ and $B_{1},\ldots,B_{m}$ be projective subspaces of an $n$-dimensional projective space $W$ over the finite  field $\mathbb{F}_{q}$. Assume
that $(A_{i},B_{i})_{1\leq i\leq m}$ is cross-intersecting. Then $$m\leq \frac{q^{n+1}-1}{q-1}-1.$$
\end{theorem}

\begin{proof} Let $t=\frac{q^{n+1}-1}{q-1}$.
For each projective subspace  $F\leq W$, its characteristic vector $v_{F}\in \{0,1\}^{t}$
 satisfies $v_{F}(s)=1$ if and only if the $s$-th projective point belongs to $F$. Here $v_{F}(s)$ denotes the $s$-th coordinate of the vector $v_{F}$.

For each $1\leq j\leq m$, let $v_{j}=(v_{j}(1),\ldots,v_{j}(t))$ and  $w_{j}=(w_{j}(1),\ldots,w_{j}(t))$ denote the characteristic vector of  $A_{j}$ and $B_{j}$ respectively.
Let $w_{m+1}=(1,1,\ldots,1)$ and $w_{m+2}=(0,0,\ldots,0)$ denote the characteristic vector of $W$ and $\emptyset$ respectively. Here $v_{j}(i)$ denotes the $i$-th coordinate of the vector $v_{j}$. Similarly, $w_{j}(i)$ denotes the $i$-th coordinate of the vector $w_{j}$.

Consider the polynomials $$P_{i}(x_{1},\ldots,x_{t}):=1-\sum_{k=1}^{t}v_{i}(k)x_{k}\in \mathbb{F}_{2}[x_{1},\ldots,x_{t}]$$ for each $1\leq i\leq m,$
  
$$P_{m+1}(x_{1},\ldots,x_{t}):=\sum_{k=1}^{t}x_{k}\in \mathbb{F}_{2}[x_{1},\ldots,x_{t}]$$

$$P_{m+2}(x_{1},\ldots,x_{t}):=\sum_{k=1}^{t}(1-x_{k})\in \mathbb{F}_{2}[x_{1},\ldots,x_{t}]$$

We claim that the polynomials $\{P_{i}\colon 1\leq i\leq m+2\}$ are linearly independent polynomial over $\mathbb{F}_{q}$. Since  $(A_{i},B_{i})_{1\leq i\leq m}$ is a cross–intersecting pair of families of projective
subspaces and the intersection of two projective subspaces is either empty or a projective subspace, we have
\begin{equation}\label{E1}
P_{i}(w_{i})=1-\sum_{k=1}^{t}v_{i}(k)w_{i}(k)=1-{\scriptsize \begin{bmatrix}
		A_{i}\cap B_{i}\\
		1 \\
	\end{bmatrix}}_q=1
\end{equation}
\begin{equation}\label{E2}
P_{i}(w_{j})=1-\sum_{k=1}^{t}v_{i}(k)w_{j}(k)=1-{\scriptsize \begin{bmatrix}
		A_{i}\cap B_{j}\\
		1 \\
	\end{bmatrix}}_q=1-\frac{q^{\ell}-1}{q-1}\equiv 0\pmod{q}
\end{equation}
\begin{equation}\label{E3}
P_{i}(w_{m+1})=1-\sum_{k=1}^{t}v_{i}(k)=1-{\scriptsize \begin{bmatrix}
		A_{i}\\
		1 \\
	\end{bmatrix}}_q=1-\frac{q^{\ell'}-1}{q-1}\equiv 0\pmod{q}
\end{equation}
\begin{equation}\label{Extra}
P_{m+2}(w_{i})=\sum_{k=1}^{t}(1-w_{i}(k))=t-\frac{q^{\ell'}-1}{q-1}\equiv 0\pmod{q}
\end{equation}
where $1\leq i\leq m$, $\ell,\ell'\geq 1$, and ${\scriptsize \begin{bmatrix}
		M\\
		1 \\
	\end{bmatrix}}_q$ counts the number of projective points in $M$. 

Consider a linear combination
\begin{equation}\label{E4}
\sum_{i=1}^{m+2}\lambda_{i}P_{i}=0,
\end{equation}
where $\lambda_{i}\in \mathbb{F}_{q}$ for $1\leq i\leq m+2$. 

 Substitute $w_{m+1}=(1,1,\ldots,1)$ for the variable of each
side of $(\ref{E4})$. Then by equation $(\ref{E3})$, all but the $(m+1)$-th  term vanish, and what
remains is $\lambda_{m+1}P_{m+1}(w_{m+1})\equiv 0 \pmod{q}.$
But $P_{m+1}(w_{m+1})=1\neq 0$ implies that $\lambda_{m+1}=0.$

Suppose $i_{0}$ is the smallest $i$ such that $\lambda_{i}\neq 0$ and $1\leq i\leq m$. Substitute $w_{i_{0}}$ for the variable of each
side of  $(\ref{E4})$. Then by equations  $(\ref{E1})$, $(\ref{E2})$ and  $(\ref{Extra})$, all but the $i_{0}$-th  term vanish, and what
remains is $\lambda_{s_{0}}P_{s_{0}}(w_{s_{0}})\equiv 0 \pmod{q}.$
But $P_{s_{0}}(w_{s_{0}})=1\neq 0$ implies that $\lambda_{s_{0}}=0$. It follows that $\lambda_{i}=0$ for any $1\leq i\leq m$.

 Substitute $w_{m+2}=(0,0,\ldots,0)$ for the variable of each
side of  $(\ref{E4})$. Then by $\lambda_{i}=0$ for all $1\leq i\leq m+1$, all but the $(m+2)$-th  term vanish, and what
remains is $\lambda_{m+2}P_{m+2}(w_{m+2})=0.$
But $P_{m+2}(w_{m+2})=1\neq 0$ implies that $\lambda_{m+2}=0$.

Hence, $P_{1},\ldots,P_{m+2}$ are linearly independent polynomials over $\mathbb{F}_{q}$. 

It is easy to see that he linearly independent polynomials $P_{1},\ldots,P_{m+2}$ are in the $\mathbb{F}_{2}$-space spanned by the monomials
$$1,x_1,\ldots,x_t.$$
Hence, $$m+2\leq t+1,$$ 
which implies that $$m\leq t-1=\frac{q^{n+1}-1}{q-1}-1.$$
\end{proof}	

\begin{remark}
Let $A_{1},\ldots,A_{m}$ and $B_{1},\ldots,B_{m}$ be projective subspaces of an $n$-dimensional projective space $W$ over finite field $\mathbb{F}_{2}$. Assume
that $(A_{i},B_{i})_{1\leq i\leq m}$ is cross-intersecting $($i.e. $A_{i}\cap B_{i}=\emptyset$ for each $1\leq i\leq m$ and $A_{i}\cap B_{j}\neq \emptyset$ whenever $1\leq i<j\leq m)$. By Theorem \ref{Main2}, we have
$$m\leq 2^{n+1}-2,$$ which implies that  Conjecture \ref{HegedusConj} of Heged{\"u}s' holds for $q=2$.
\end{remark}

\section{conclusion}
In this note we have provided an precise bound for cross-interecting affine subspaces, we have also proved Heged{\"u}s‘ conjecture holds for $q=2$. The conjecture of cross-intersecting projective subspaces is still open for $q\neq 2$. It is meaningful to prove that the conjecture holds for any field, and to find a concrete construction of cross-intersecting projective subspaces with size $2^{n+1}-2$. 

\bigskip
\noindent {\bf Data availability}

No data was used for the research described in this article.
	
\bigskip
\noindent {\bf Conflict of Interest}

The authors declare that they have no conflicts of interest.

\end{document}